\documentclass[11pt,a4paper]{article}
\usepackage[T1]{fontenc}
\usepackage{lmodern,amsmath,amssymb,amsthm,graphicx,booktabs,microtype}
\usepackage[margin=24mm]{geometry}
\usepackage{float,authblk,needspace}
\usepackage[hidelinks]{hyperref}
\graphicspath{{output/pdf/task5_figures/}{output/pdf/task6_figures/}{../output/pdf/task5_figures/}{../output/pdf/task6_figures/}}

\newcommand{\D}{\mathcal D}
\newcommand{\F}{\mathcal F}
\newcommand{\AND}{\mathop{\bigwedge_R}\nolimits}
\newcommand{\OR}{\mathop{\bigvee_R}\nolimits}
\DeclareMathOperator{\dist}{dist}

\newtheorem{proposition}{Proposition}
\title{R-Function Representation and Composition\\for Flexibility Analysis and Robust Design Spaces}
\author[1]{S. Kucherenko}
\author[2]{S. Pistikopoulos}
\author[3]{L. Lamakani}
\affil[1]{Imperial College London, SW7 2AZ, UK}
\affil[2,3]{Texas A\&M University, TX 77843-3251, USA}

\date{}
\begin{document}
\maketitle
\vspace{-1.5em}
\begin{abstract}
Explicit multiparametric solutions provide analytical feasibility and design-space descriptions after elimination of adjustable variables. R-functions convert these descriptions into fields whose nonnegative regions represent feasible or robust designs. The principal contribution is the analytical composition of independently derived design spaces through union, intersection with engineering restrictions, and exclusion of prohibited regions, without repeating the underlying multiparametric elimination. A support-function construction supplies the robust design fields used in this composition. A two-parameter example demonstrates how alternative configurations, a manufacturing restriction and an exclusion combine into one admissible field. The resulting representations also support geometric interrogation of external boundaries, internal interfaces and limiting uncertainty contact, illustrated by a four-parameter heat-exchanger network. Explicit MPC offers a potential further application.
\end{abstract}
\section{Introduction}
Multiparametric programming can eliminate adjustable operating variables and provide explicit analytical descriptions of feasible and design spaces~\cite{mpcontrol}. In flexibility analysis, Bansal et al.~\cite{bansal} obtain feasibility functions, flexibility indices and robust design spaces for linear systems. These projected descriptions retain the ability of operating variables to adjust to each uncertainty realization while the design remains fixed.

Once design spaces have been derived independently for alternative process configurations or technologies, subsequent decisions may require their union, intersection with additional engineering restrictions, or removal of prohibited regions. A conventional conjunction of inequalities naturally represents an intersection of feasible half-spaces or constraint sets. R-functions provide a single analytical field for unions, intersections and exclusions of independently derived design spaces. This composition does not require introducing discrete selector variables into the already-derived design-space representation.

R-function theory~\cite{rvachev,shapiro} supplies sign-preserving functional Boolean operations, with a recent application to feasible-region identification by Kucherenko et al.~\cite{rfeas}. Here these operations connect explicit multiparametric solutions to reusable analytical design-space fields. The central step is composition: once component fields have been obtained, alternatives, restrictions and exclusions can be combined without repeating unchanged multiparametric elimination.

The contribution has three parts. First, we construct analytical fields $R_F(\theta,d)$ and $R_D(d;t)$ for projected feasible and robust design spaces. Second, we compose independently derived robust design spaces through union, intersection/restriction and exclusion, with configuration-selection and uncertainty quantifiers made explicit. Third, we interpret the resulting geometry through external boundaries, internal interfaces, Euclidean clearance and higher-dimensional limiting uncertainty contact. Section~\ref{sec:composition} presents the central composition result and its concrete demonstration; Section~\ref{sec:mpc} sketches a possible explicit-MPC application.

\Needspace{5\baselineskip}
\section{R-functions after parametric feasibility elimination}
Let $d$ denote fixed design variables, $\theta$ uncertain parameters, and $y$ the state and adjustable operating variables. The projected feasible set is
\begin{equation}
 \F(d)=\{\theta:\exists y,\ h(y,\theta,d)=0,\ g_j(y,\theta,d)\le0\ (j=1,\ldots,m)\}.
 \label{eq:projection}
\end{equation}
Thus $\F(d)$ contains all uncertainty realizations $\theta$ for which at least one feasible adjustable operating vector $y$ exists. Variable bounds are included in these conditions. Suppose parametric optimization supplies a finite, attained feasibility value
\begin{equation}
 \psi(\theta,d)=\min_{y:h=0}\max_j g_j(y,\theta,d)
                  =\max_{k=1,\ldots,K}\psi^k(\theta,d)
 \label{eq:envelope}
\end{equation}
throughout the domain under consideration. The last equality is a global envelope assumption, satisfied by the illustrative linear model below. Set $\phi_k=-\psi^k$. Write $\AND$ for repeated R-conjunction, an operation representing intersection through nonnegative primitive fields (specific choices follow below). Then
\begin{equation}
 R_F(\theta,d)=\AND_{k=1}^K\phi_k(\theta,d),\qquad
 \theta\in\F(d)\ \Longleftrightarrow\ R_F(\theta,d)\ge0.
 \label{eq:RF}
\end{equation}
For the minimum R-conjunction, the value identity is exact:
\[
 R_F=\min_k\phi_k=-\max_k\psi^k=-\psi.
\]
For a different sign-equivalent conjunction, such as $\wedge_0$ below,
$R_F\ge0\iff\psi\le0$, but generally $R_F\ne-\psi$. Parametric programming supplies the projected feasibility description after elimination of the adjustable variables; the R-function representation then makes this description available for subsequent analytical composition.

For positive-inside primitives, two standard choices are~\cite{rvachev,shapiro}
\begin{align}
 a\wedge_1 b&=\min(a,b),& a\vee_1 b&=\max(a,b),\label{eq:minmax}\\
 a\wedge_0 b&=a+b-\sqrt{a^2+b^2},&
 a\vee_0 b&=a+b+\sqrt{a^2+b^2}.\label{eq:algebraic}
\end{align}
Negation is $\neg_R a=-a$. Both conjunctions are nonnegative exactly when both arguments are nonnegative; both disjunctions are nonnegative exactly when at least one argument is nonnegative. The operator $\OR$ denotes repeated R-disjunction. Multiargument operations are left-associated in the displayed order. They are exact sign constructions, not approximations to minimum or maximum.

For locally valid branches, fix $d$ and define an analysis domain
\[
 \mathcal A(d)=\bigcup_{k=1}^K C_k(d),\qquad
 C_k(d)=\{\theta:R_{C_k}(\theta,d)\ge0\}.
\]
Assume that the regions cover the intended domain and that, on every region,
\[
 \theta\in C_k(d):\quad
 \theta\in\F(d)\ \Longleftrightarrow\ \psi^k(\theta,d)\le0.
\]
A valid global representation of $\F(d)\cap\mathcal A(d)$ is then
\begin{equation}
 R_F(\theta,d)=\OR_k\bigl[R_{C_k}(\theta,d)\wedge_R(-\psi^k(\theta,d))\bigr].
 \label{eq:regional}
\end{equation}
Every local branch is clipped to its region. Values used outside that region only permit evaluation of the expression and have no feasibility meaning. On overlaps and common boundaries, the assumed equivalence requires consistent feasibility classification; duplicate coverage is harmless for the represented nonnegative set under R-disjunction, although the numerical field value may depend on the representation and composition. For closed critical regions, shared boundaries are included. The intersection shortcut~\eqref{eq:RF} requires globally valid branches and must not be inferred merely from local affine formulas.

The regional union~\eqref{eq:regional} may have zeros on internal interfaces. Exact membership does not imply either strict interior positivity or a distance interpretation.

\Needspace{5\baselineskip}
\section{Flexibility as uncertainty-set containment}
Let $T(\delta)=\bar\theta+\delta U$, where $U$ is compact and convex, $0\in U$, and $\delta\ge0$. Thus $0\le\delta_1\le\delta_2$ implies $T(\delta_1)\subseteq T(\delta_2)$. On a domain where~\eqref{eq:RF} is valid, flexibility is equivalently
\begin{align}
 T(\delta)&\subseteq\{\theta:R_F(\theta,d)\ge0\},\label{eq:contain}\\
 F(d)&=\sup\{\delta\ge0:R_F(\theta,d)\ge0\quad\forall\theta\in T(\delta)\}.
 \label{eq:F}
\end{align}
Here $F(d)$ denotes the scalar flexibility index and should be distinguished from the projected feasible set $\F(d)$ defined in~\eqref{eq:projection}.
The containment $T(\delta)\subseteq\F(d)$ means that for every $\theta\in T(\delta)$ there exists a feasible $y=y(\theta)$. These adjustable operating variables are recourse chosen after the uncertainty realization is known; the design stays fixed. It does not require one operating decision to work for all uncertainty realizations. We assume nominal feasibility when using~\eqref{eq:F}; an infeasible nominal point gives an empty admissible-scale set and requires an explicit reporting convention.

The field $R_F$ gives one common functional description independent of the uncertainty shape subsequently tested against it. Its zero set provides candidate boundary/contact geometry; internal partition-interface zeros must be distinguished from the external boundary. R-conjunction and R-disjunction also permit nonconvex or disconnected feasible sets assembled from multiple components. Thus the same representation can be composed with additional restrictions and tested against different uncertainty sets.

For globally valid affine branches, support-function robustification reduces containment to explicit inequalities, as developed in the next section. R-functions retain the resulting robust design geometry as reusable fields for subsequent analytical composition.

For this uncertainty shape, define $h_U(a)=\sup_{u\in U}a^Tu$. For a symmetric box with half-widths $\Delta_j$, $h_U(a)=\sum_j|a_j|\Delta_j$; for the unit Euclidean ball, $h_U(a)=\|a\|_2$. These different supports yield different flexibility margins for the same feasible geometry. If the parametric model is valid only within a hard parameter window, admissible uncertainty sets must also remain inside that window.

\Needspace{4\baselineskip}
\subsection*{A concise linear illustration}
For the linear system of Bansal et al.~\cite{bansal}, with
$\theta\in[0,4]^2$ and $d\in[0,5]^2$, the established feasibility envelope has two branches:
\begin{align}
 \psi^1&=-\tfrac23\theta_1-\tfrac14\theta_2+\tfrac12d_1-d_2+\tfrac16,\notag\\
 \psi^2&=\tfrac13\theta_1-\theta_2-\tfrac12d_1+\tfrac12d_2-\tfrac13.
 \label{eq:example}
\end{align}
The corresponding field is $R_F=(-\psi^1)\wedge_R(-\psi^2)$.
Take $T(\delta)=[2-2\delta,2+2\delta]^2$, with $0\le\delta\le1$ so that the stated window is respected. At $d=(3,1)$, the known index is $F=7/11$~\cite{bansal}. The limiting box touches $R_F=0$ at $(8/11,8/11)$, where $\psi^1=0$ and $\psi^2=-20/11$. The contact point gives a geometric interpretation of this flexibility index.

For this example, the two affine supports produce the design inequalities developed next. Containment has thereby converted a field in uncertainty-design space into a field in design space indexed by the requested uncertainty scale. This is the main transition from pointwise feasible geometry to robust design geometry.

\Needspace{5\baselineskip}
\section{R-function representation of robust design spaces}
Let $\mathcal B$ be the admissible design domain, with exact field $R_B$. The central object is $\D(t)=\{d\in\mathcal B:T(t)\subseteq\F(d)\}$: each design certifies feasible adjustable operation throughout the uncertainty set.
\begin{proposition}[Robust design field after affine elimination]
Suppose projected feasibility is equivalent to the globally valid inequalities $a_k^T\theta+c_k(d)\le0$, $k=1,\ldots,m$, throughout the relevant domain. Let $T(t)=\bar\theta+tU$, where $t\ge0$ and $U$ is compact and convex and contains the origin. Define
\begin{equation}
 q_k(d;t)=-[a_k^T\bar\theta+c_k(d)+t h_U(a_k)].
 \label{eq:support}
\end{equation}
For targets whose uncertainty sets remain within the validity domain,
\begin{equation}
 \D(t)=\{d\in\mathcal B:q_k(d;t)\ge0\ \forall k\}
       =\{d:R_D(d;t)\ge0\},\quad
 R_D=R_B\wedge_R\left(\AND_k q_k\right).
 \label{eq:RD}
\end{equation}
If no separate design restrictions apply, the $R_B$ factor is omitted.
\end{proposition}
\begin{proof}
For each $k$, compactness ensures that
$\sup_{\theta\in T(t)}(a_k^T\theta+c_k(d))
=a_k^T\bar\theta+c_k(d)+t h_U(a_k)$.
Universal feasibility is therefore equivalent to every $q_k\ge0$. R-conjunction gives exactly their intersection with $\mathcal B$. The existential recourse quantifier does not appear explicitly in~\eqref{eq:support}--\eqref{eq:RD}, because the adjustable variables have already been eliminated in deriving the projected feasibility inequalities.
\end{proof}
The resulting field $R_D(d;t)$ provides an analytical representation of the complete family of robust design spaces and forms the basic object for the composition operations in Section~\ref{sec:composition}. Since $0\in U$, $h_U(a_k)\ge0$. Every robust inequality therefore becomes monotonically more restrictive as $t$ increases: $t_2\ge t_1\ge0$ implies $\D(t_2)\subseteq\D(t_1)$ for admissible targets.

For the established example of Bansal et al.~\cite{bansal}, $\mathcal B=[0,5]^2$ and $0\le t\le1$. Positive rescaling of the general primitives gives
\begin{align}
 Q_1&=-3d_1+6d_2+10-11t,&Q_2&=3d_1-3d_2+10-16t,\notag\\
 R_D&=\AND\left(d_1,5-d_1,d_2,5-d_2,Q_1/\sqrt{45},Q_2/\sqrt{18}\right).
 \label{eq:q}
\end{align}
Here the repeated R-conjunction $\AND$ combines six normalized primitive fields: $R_D\ge0$ is equivalent to $d_1\ge0$, $5-d_1\ge0$, $d_2\ge0$, $5-d_2\ge0$, $Q_1\ge0$ and $Q_2\ge0$. At fixed $t$, for an affine field $Q(d)=a^Td+c$ with $a\ne0$, $Q/\|a\|_2$ is signed Euclidean distance to the hyperplane $Q=0$, positive on $Q>0$. The norms $\|(-3,6)\|_2=\sqrt{45}$ and $\|(3,-3)\|_2=\sqrt{18}$ therefore normalize $Q_1,Q_2$ consistently with the four unit-normal box primitives.

At $t=1$, the active boundaries are $Q_1=0$, $Q_2=0$ and $d_1=5$. Their pairwise intersections give $(13/3,7/3)$, $(5,8/3)$ and $(5,3)$, respectively forming the established robust triangle; its first vertex minimizes $10(d_1+d_2)$~\cite{bansal}. This robust design field supplies a component for the composition developed next.
\begin{figure}[!htbp]\centering
\includegraphics[width=.90\linewidth]{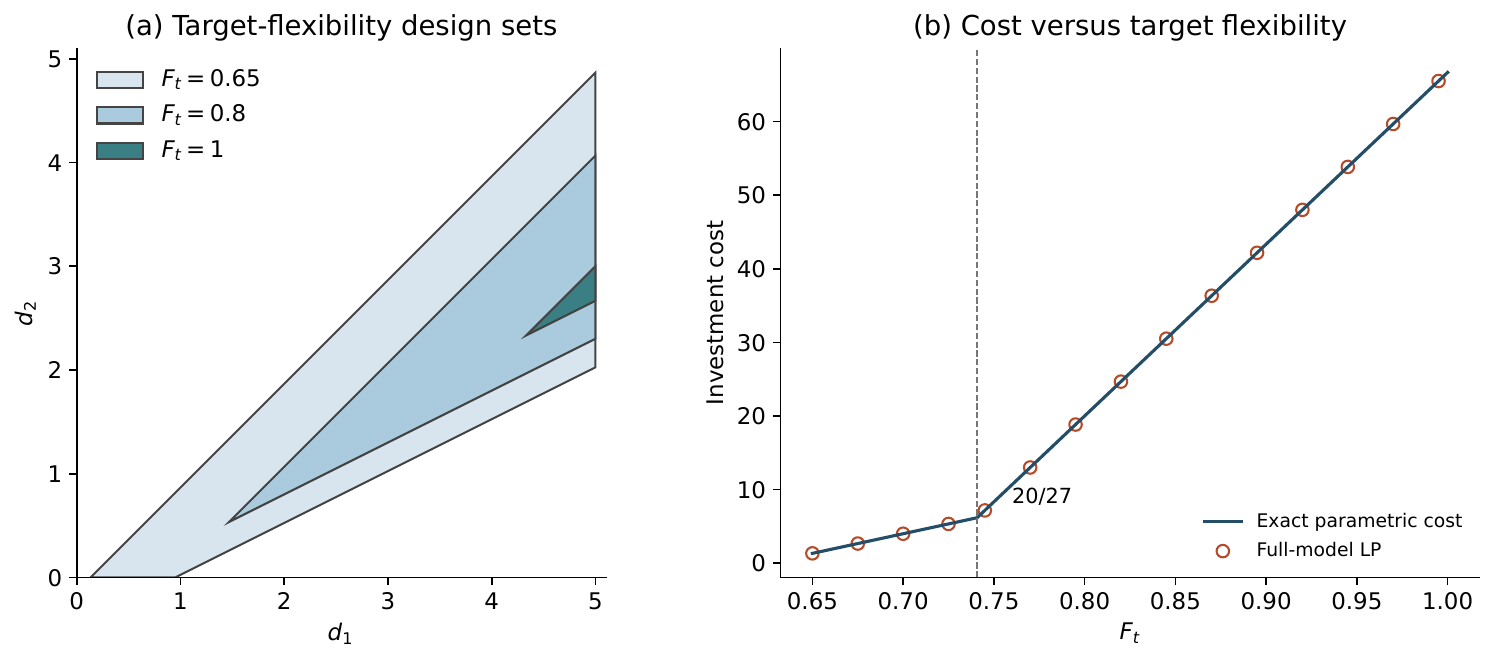}
\caption{The robust design family and established flexibility--cost relationship of Bansal et al.~\cite{bansal}; $F_t$ denotes target $t$. Circles denote independent LP solutions. The R-field represents the same feasible choices, each certified over a complete uncertainty set.}
\label{fig:targets}
\end{figure}

\Needspace{5\baselineskip}
\section{Analytical composition of robust design spaces}\label{sec:composition}
Suppose separate multiparametric analyses have produced robust design fields $R_{D_1}(d;t)$ and $R_{D_2}(d;t)$ for two alternative configurations or technologies. Once available in common design coordinates, these fields can be combined analytically without rerunning the underlying elimination. If $R_S$ describes an independent safety or design restriction, their composition is
\begin{equation}
 R_A(d;t)=\bigl(R_{D_1}(d;t)\vee_R R_{D_2}(d;t)\bigr)\wedge_R R_S(d),
 \label{eq:alternative}
\end{equation}
which represents $(\D_1(t)\cup\D_2(t))\cap\D_S$ exactly. Each alternative is analysed separately. The designer selects an alternative $i$ for which feasible operating variables $y_i(\theta)$ exist for every uncertainty realization in $T(t)$. The operating variables may change with $\theta$; the alternative itself is not switched after $\theta$ is observed. Thus the union permits either alternative, provided that the selected one is robust for the complete uncertainty set. This differs from allowing selection after observation:
\[
 (\exists i\ \forall\theta\ \exists y_i)\quad\ne\quad
 (\forall\theta\ \exists i\ \exists y_i).
\]
The R-union is therefore applied to the robust design fields, after each alternative has been certified over the complete uncertainty set.

Let an exclusion field $R_E$ be positive precisely on the forbidden open region $\D_E$. Then
\begin{equation}
 R_{\rm adm}(d;t)=R_A(d;t)\wedge_R(-R_E(d))
 \label{eq:exclusion}
\end{equation}
represents $(\D_1(t)\cup\D_2(t))\cap\D_S\setminus\D_E$ using
$R_{\rm adm}\ge0$. The boundary $R_E=0$ is retained. If the forbidden region is instead the closed set $\{R_E\ge0\}$, exact subtraction requires $R_E<0$ together with $R_A\ge0$. A continuous nonnegative-superlevel representation cannot in general encode that nonclosed difference. This convention matters for contact with exclusion boundaries. Together,~\eqref{eq:alternative} and~\eqref{eq:exclusion} provide a single analytical representation of a composite design space assembled from independently obtained multiparametric spaces. Additional restrictions and exclusions act directly on these reusable fields.

\Needspace{4\baselineskip}
\subsection*{Illustrative configuration and design restrictions}
\begin{figure}[!htbp]\centering
\includegraphics[width=.80\linewidth]{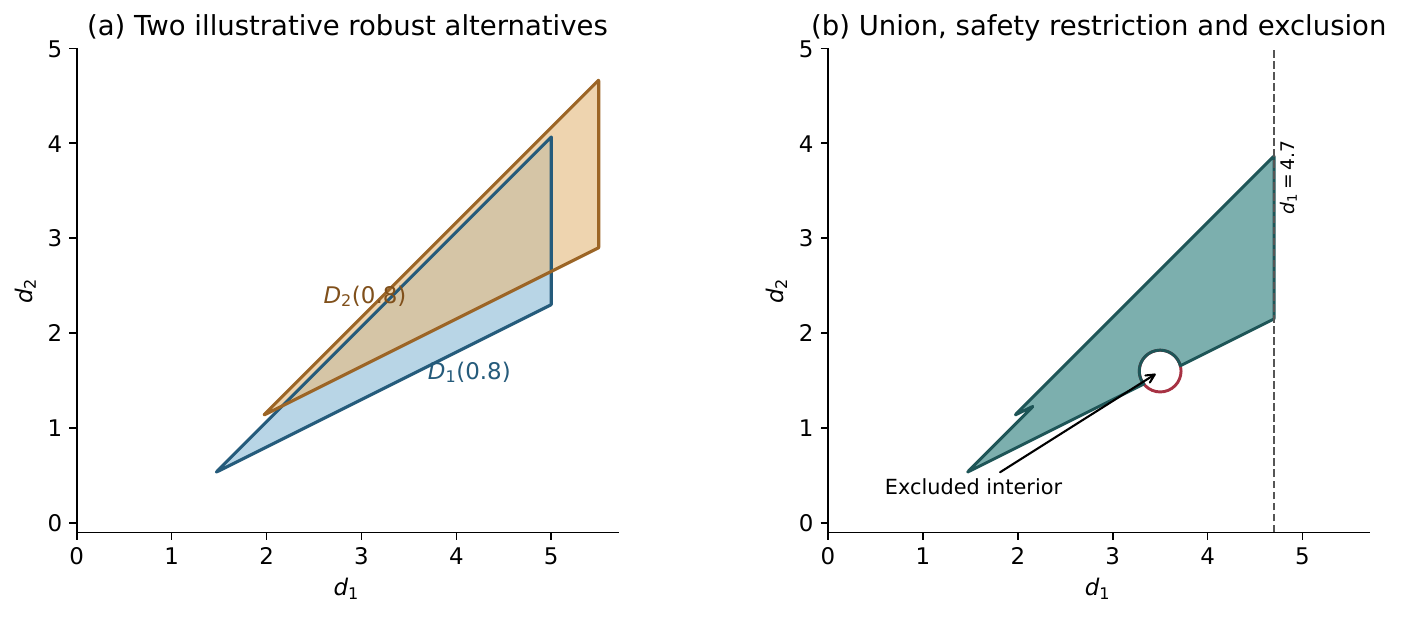}
\caption{Analytical composition at target $t=0.8$: (a) component robust design spaces; (b) their R-disjunction union, intersected with a manufacturing restriction and modified by exclusion of an open disk. The final admissible region is represented by one field $R_{\rm adm}(d;t)$. The disk circumference remains admissible wherever the other constraints hold.}
\label{fig:compose}
\end{figure}

Interpret $\D_1(t)$ and $\D_2(t)$ as robust specifications for technologies A and B, $\D_S$ as an independent manufacturing limit, and $\D_E$ as a prohibited design combination. This is an illustrative construction, not a calibrated engineering example. Take
$\D_1(t)=\D(t)$, $\D_2(t)=\D(t)+(0.5,0.6)$ and
$R_{D_2}(d;t)=R_D(d-(0.5,0.6);t)$. The translation specifies an offset in design coordinates. Use $R_S=4.7-d_1$ and the prohibited open disk
\[
 R_E(d)=0.22^2-(d_1-3.5)^2-(d_2-1.6)^2>0.
\]
At $t=0.8$,~\eqref{eq:exclusion} yields the nonconvex admissible set in Figure~\ref{fig:compose}. The sequence mirrors the reuse of independently derived spaces: R-disjunction forms the union, R-conjunction imposes the manufacturing limit, and negation of the exclusion field removes the prohibited disk. The final region remains available as the single analytical field $R_{\rm adm}(d;t)$.

\Needspace{5\baselineskip}
\section{Geometrical margins and boundary interpretation}\label{sec:margins}
Analytical fields produced by composition support both membership queries and geometric interpretation. For convex components, normalization relates primitive slack to Euclidean clearance; for composite regions, the decomposition also identifies internal interfaces.
Let $D=\{d:a_j^Td\le b_j,\ j=1,\ldots,m\}$ be a nonempty, full-dimensional, proper closed convex polyhedron, with nonzero normals. Define
\begin{equation}
 \rho_j(d)=\frac{b_j-a_j^Td}{\|a_j\|_2},\qquad
 r_D(d)=\min_j\rho_j(d).
 \label{eq:rho}
\end{equation}
For $d\in D$,
\begin{equation}
 r_D(d)=\dist_2(d,\partial D).
 \label{eq:distance}
\end{equation}
Indeed, a Euclidean ball of radius $r$ centered at $d$ lies in $D$ exactly when $a_j^Td+r\|a_j\|_2\le b_j$ for every $j$. The largest such radius is both the minimum normalized slack and the distance to the boundary. Individual $\rho_j$ measure distance to supporting hyperplanes; their orthogonal projections need not lie on the corresponding finite facets. The ball argument establishes~\eqref{eq:distance} nevertheless, including redundant valid inequalities.

For nonnegative $a,b$, $0\le a\wedge_0b\le\min(a,b)$. Thus the left-associated algebraic R-conjunction satisfies $0\le R_D^0(d)\le\dist_2(d,\partial D)$ for $d\in D$. This is a lower bound on interior Euclidean clearance for the normalized convex representation, not itself a Euclidean distance or radius. At the incenter,
\begin{align*}
 d_c&\approx(4.890,2.735),\\
 \dist_2(d_c,\partial\D(1))&\approx0.110,
 &R_D^0(d_c;1)&\approx0.046.
\end{align*}
The two fields represent the same set. Define signed distance, positive inside, by
\[
 s_D(x)=\begin{cases}
 \dist_2(x,\partial D),&x\in D,\\
 -\dist_2(x,D),&x\notin D.
 \end{cases}
\]
The notation $s_D$ is distinct from the design point $d_s$ below. For $D=\D(1)$ and $x=(5.2,3.2)$, $r_D(x)=-0.2$ but $s_D(x)=-\sqrt{0.08}\approx-0.283$: the normalized minimum is not a global signed distance.

Design clearance $\dist_2(d,\partial\D(t))$ measures design-variable perturbations at fixed $t$; flexibility $F(d)$ measures prescribed uncertainty expansion in parameter space. They are not interchangeable. At $(13/3,7/3)$, $F=1$ but clearance to $\partial\D(1)$ is zero. Thus a design can tolerate the full uncertainty box without positive isotropic design clearance, which also depends on coordinate scaling.

\Needspace{4\baselineskip}
\subsection*{External boundaries versus internal partition interfaces}
To show why a composite R-field zero need not be an external-boundary point, artificially partition the robust triangle without changing its overall set. The auxiliary field $\sigma(d)=(81d_1-129d_2-50)/\sqrt{81^2+129^2}$ is signed perpendicular distance to the partition line $\sigma=0$; its positive and negative signs specify the two sides. The partition isolates the effect of an internal interface in a composed field. With $r=R_D^{\min}(d;1)$, the union field
\begin{equation}
 R_{\rm split}(d)=\max\{\min(r,\sigma),\min(r,-\sigma)\}
 \label{eq:interface}
\end{equation}
represents the union of the two subregions, exactly the original triangle, yet vanishes on their artificial internal interface. At the exact interior point $d_s=(24/5,1694/645)\approx(4.800,2.626)$, the partition field is zero, whereas the external-boundary clearance is $\dist_2(d_s,\partial\D(1))\approx0.053$. The internal zero does not indicate loss of robust feasibility.

For a general R-expression,
\[
 \boxed{R(x)=0\ \not\Rightarrow\ x\in\partial\{R\ge0\}.}
\]
Boundary interpretation depends on the functional decomposition. External clearance excludes internal interface zeros; neither a metric interpretation nor the convex minimum identity transfers automatically to arbitrary unions or R-expressions.

\Needspace{5\baselineskip}
\section{Higher-dimensional illustration}
\Needspace{4\baselineskip}
\subsection*{Four uncertain temperatures}
The established two-hot/two-cold heat-exchanger network of Bansal et al.~\cite{bansal} illustrates the geometry and dimension of limiting uncertainty contact at its flexibility index $F=1/2$. The uncertain temperatures are $\theta=(T_1,T_3,T_5,T_8)$. Their nominal values are $(620,388,583,313)$ K, each with a deviation of $\pm10$ K. After elimination of the adjustable variables, its established feasibility function is $\psi_H=\max_k p_k$, where
\begin{align}
 p_1&=-\tfrac35T_1-\tfrac15T_3-\tfrac25T_5-\tfrac45T_8+922,\notag\\
 p_2&=\tfrac12T_8-\tfrac{323}{2},\notag\\
 p_3&=-\tfrac37T_1-\tfrac17T_3-\tfrac47T_5+\tfrac{4504}{7},\notag\\
 p_4&=-\tfrac13T_5+T_8-\tfrac{376}{3}.
 \label{eq:hen}
\end{align}
Consequently $R_H=\AND_{k=1}^4(-p_k)$ represents the projected feasible temperature set, with the analysis window understood. With no design variables, this example applies the containment construction directly in higher-dimensional uncertainty space.

For $T_H(\delta)=\bar\theta+10\delta[-1,1]^4$, the known flexibility index is $F=1/2$~\cite{bansal}. The four branch contact scales are $14/25$, $1$, $19/20$ and $1/2$. At $\delta=1/2$, the box ranges are $615\le T_1\le625$, $383\le T_3\le393$, $578\le T_5\le588$ and $308\le T_8\le318$. The limiting constraint is $p_4$, with contact at $T_5=578$, $T_8=318$. Since $p_4$ is independent of $T_1$ and $T_3$, the entire two-dimensional face
\begin{equation}
 \{(T_1,T_3,578,318):615\le T_1\le625,\ 383\le T_3\le393\}
 \label{eq:face}
\end{equation}
is contained in $T_H(1/2)\cap\{R_H=0\}$. At that scale the other three branches remain strictly feasible throughout the box, with maxima $-1.2$, $-2.5$ and $-36/7$. Within the uncertainty box $T_H(1/2)$, the limiting feasibility contact is precisely the face selected by $p_4=0$, rather than a unique critical point. A scalar flexibility index gives the limiting scale, beyond which feasibility is lost, but does not distinguish an isolated point, edge, face or higher-dimensional contact set. Here $T_H(F)\cap\{R_H=0\}$ retains that geometric information. Because all other branches are strictly feasible, these zeros are genuine feasibility-boundary contacts, not internal partition zeros. The analytical field thus retains the geometry and dimension of limiting uncertainty contact beyond the scalar flexibility index.

\Needspace{5\baselineskip}
\section{Potential Extension to Explicit MPC}\label{sec:mpc}
\begin{center}
\small
\begin{tabular}{c@{\quad$\longrightarrow$\quad}c@{\quad$\longrightarrow$\quad}c}
$\{CR_i\}$ & $R_i(x)$ & $R_{\rm MPC}(x)=\OR_i R_i(x)$\\[3pt]
Critical regions & Region fields & Feasible-domain union
\end{tabular}
\end{center}
An explicit MPC solution supplies critical regions and associated local affine control laws~\cite{mpcontrol}. The same R-function construction could combine the critical-region geometry into a single analytical representation of the feasible-domain union, making it available for subsequent geometric interrogation. Internal switching interfaces should be distinguished from the external boundary of that union, as in Section~\ref{sec:margins}. This offers a potential future application of the design-space composition framework.

\section{Conclusions}
Starting from explicit multiparametric solutions, this work constructs analytical R-function fields for projected feasible and robust design spaces. Its principal contribution is to make independently derived design spaces reusable analytical objects: their union, intersection with additional engineering restrictions, and exclusion of prohibited regions can be assembled into one admissible field without repeating unchanged multiparametric elimination. The composition example demonstrates this sequence for a nonconvex design region; the same operations accommodate disconnected component spaces.

The analytical representation also supports geometric interrogation. Normalized convex fields quantify design clearance, while the functional decomposition distinguishes external boundaries from internal interfaces. In uncertainty space, limiting-contact sets retain geometric and dimensional information beyond the scalar flexibility index. These capabilities connect multiparametric feasibility analysis to subsequent design-space composition and interpretation, with explicit MPC as a possible extension.

\section*{Acknowledgements}
S.K. gratefully acknowledges Prof. Nilay Shah of Imperial College London for his interest in and support of this work.

\end{document}